\documentclass[12pt]{article}

\usepackage{amsmath}
\usepackage{amsfonts}
\usepackage{amssymb}
\usepackage{extarrows}
\usepackage{plain}
\usepackage{amsthm}
\usepackage{mathrsfs}
\usepackage{graphicx}
\usepackage[noend]{algpseudocode}
\usepackage{multirow}
\newtheorem{theorem}{Theorem}[section]
\newtheorem{lemma}{Lemma}[section]
\newtheorem{coro}[theorem]{\bf Corollary}

\newtheorem{remark}[theorem]{\bf Remark}

\newtheorem{defn}[theorem]{\bf Definition}

\def\endproof{{\hfill $\square$}\medskip}

\newcommand{\gd}{{\bigtriangledown}}

\newcommand{\del}{\nabla}
\newcommand{\laplace}{\triangle}
\newcommand{\Cov} {{\rm Cov}}
\newcommand{\COV} {{\rm COV}}

\AtEndDocument{\bigskip{\footnotesize
        \textsc{Faculty of Innovation Engineering, Macau University of Science and Technology, Macao, China}\par
		\textit{E-mail address: }\texttt{pdang@must.edu.mo}
		\par
		\medskip
		\textsc{Macao Center for Mathematical Sciences, Macau University of Science and Technology, Macao, China}\par
		\textit{E-mail address: }\texttt{wxmai@must.edu.mo}
		
}}

\begin{document}



\title{Improved Caffarelli-Kohn-Nirenberg Inequalities and Uncertainty Principle}


\author{Pei Dang, Weixiong Mai\footnote{ Corresponding author, wxmai@must.edu.mo.}}

\maketitle	

%\address{1. Faculty of Information Technology, Macau University of Science and Technology, Macao, China\\
%2. Macao Center for Mathematical Sciences, Macau University of Science and Technology, Macao, China. Email: wxmai@must.edu.mo 
%}

\begin{abstract}
In this paper we prove some improved Caffarelli-Kohn-Nirenberg inequalities and uncertainty principle for complex- and vector-valued functions on $\mathbb R^n$, which is a further study of the results in \cite{Dang-Deng-Qian}. In particular, we introduce an analogue of \lq\lq phase derivative" for vector-valued functions. Moreover, using the introduced \lq\lq phase derivative", we extend the extra-strong uncertainty principle to cases for complex- and vector-valued functions defined on $\mathbb S^n,n\geq 2.$

\end{abstract}

	\section{Introduction}
	Denote by $C_0^\infty(\mathbb R^n)$ the smooth functions with compact support in $\mathbb R^n$.  Let $0<q<2<p \text{ and }2<n<\frac{2(p-q)}{p-2}$.
	A subclass of the well-known Caffarelli-Kohn-Nirenberg (CKN) inequalities \cite{Caffarelli-Kohn-Nirenberg} states that for any real-valued $f\in C_0^\infty(\mathbb R^n)$ there holds
	\begin{align}\label{CKN-RN}
		\int_{\mathbb R^n}|\del f(x)|^2dx \int_{\mathbb R^n} \frac{|f(x)|^{2p-2}}{|x|^{2q-2}}dx\geq \frac{(n-q)^2}{p^2}\int_{\mathbb R^n}\frac{|f(x)|^p}{|x|^{q}}dx.
	\end{align}
	The constant $\frac{(n-q)^2}{p^2}$ is sharp (see \cite{Xia}). See also \cite{Nguyen} for a more general discussion on the sharp constant of CKN inequalities. The general CKN inequalities contain many famous inequalities, for examples, Sobolev inequalities, Hardy inequalities, the Gagliardo-Nirenberg inequalities and the Heisenberg-Pauli-Weyl (HPW) uncertainty principle (see e.g. \cite{Kombe-Ozaydin,Kombe-Ozaydin1,Hei,Xia,Xia1,Xia2}). Therefore, CKN inequalities play an important role in the theory of partial differential equations. Moreover, finding sharp constants and extremals for CKN inequalities on Riemannian manifolds are interesting and non-trivial problems (see e.g. \cite{Aubin,Kri,Kri1,Kri-Ohta,Xia,Nguyen}).
	
Since, as mentioned above, the HPW uncertainty principle can be viewed as a special case of CKN inequalities and HPW uncertainty principle has undergone an important development in different settings in recent years, that is essentially represented by formula (\ref{ex-strong}), we wonder whether there exist analogous results for CKN inequalities. It is the problem that motivates us  to consider CKN inequalities, and, furthermore, obtain the results in the present paper. Now we first review recent developments of HPW uncertainty principle. HPW uncertainty principle plays an important role in quantum mechanics and harmonic analysis, which states that the position and the momentum of a particle cannot be both determined precisely in any quantum state (see e.g. \cite{Hei,Weyl,Fefferman,Folland-Sitaram}). Or equivalently, in the language of harmonic analysis, the HPW uncertainty principle says that a nonzero function and its Fourier transform cannot both be sharply localized. In \cite{Gabor} Gabor brings uncertainty principle to the sight of signal analysts. From the aspect of signal analysis, the HPW uncertainty principle has also been extensively studied, and there are some new developments in recent years.
	It is stated as follows,
	\begin{align}\label{UP-classic}
		\sigma_{t,s}^2\sigma_{\omega,s}^2\geq \frac{1}{4},
	\end{align}
	where $\sigma_{t,s}$ and $\sigma_{\omega,s}$, respectively, denote the duration and the bandwidth of a signal $s(t)\in L^2(\mathbb R)$ with $||s||_{L^2}=1$, and are defined as
	\begin{align*}
		\sigma_{t,s}^2 = \int_{-\infty}^\infty (t-\langle t \rangle_s)^2 |s(t)|^2 dt
	\end{align*}and
	\begin{align*}
		\sigma_{\omega,s}^2 = \int_{-\infty}^\infty (\omega-\langle \omega\rangle_s)^2|\hat s(\omega)|^2 d\omega,
	\end{align*}
	where $\langle t\rangle_s$ and $\langle \omega\rangle_s $ are the means of time $t$ and Fourier frequency $\omega$, given by
	$\langle t\rangle_s = \int_{-\infty}^\infty t|s(t)|^2dt$ and $\langle \omega\rangle_s=\int_{-\infty}^\infty \omega |\hat s(\omega)|^2 d\omega,$ respectively. Here $\hat s(\omega)$ denotes the Fourier transform of $s$, i.e., if $s(t) \in L^1(\mathbb R)$, then $\hat s(\omega)= \frac{1}{\sqrt{2\pi}}\int_{-\infty}^\infty s(t)e^{-i t\omega}dt.$ For a signal written as $s(t)=|s(t)|e^{i\varphi(t)},$ stronger versions of uncertainty principle are available (cf. \cite{Cohen, Dang-Deng-Qian}), which are stated as follows,
	\begin{align}\label{strong}
		\sigma_{t,s}^2\sigma_{\omega,s}^2\geq \frac{1}{4}+ \Cov_s^2
	\end{align}
	and
	\begin{align}\label{ex-strong}
		\sigma_{t,s}^2\sigma_{\omega,s}^2\geq \frac{1}{4} +\COV_s^2,
	\end{align}
	where
	\begin{align*}
		\Cov_s = \int_{-\infty}^\infty (t-\langle t\rangle_s)(\varphi^\prime(t)-\langle \omega\rangle_s)|s(t)|^2 dt
	\end{align*}
	and
	\begin{align*}
		\COV_s = \int_{-\infty}^\infty |(t-\langle t\rangle_s)(\varphi^\prime(t)-\langle \omega\rangle_s)| |s(t)|^2 dt.
	\end{align*}
	It is obvious that $\COV^2_s\geq \Cov_s^2.$ In the literature (\ref{strong}) is called \lq\lq strong uncertainty principle", and (\ref{ex-strong}) is called \lq\lq extra-strong uncertainty principle".
	Recently, the study of (\ref{UP-classic}), (\ref{strong}) and (\ref{ex-strong}) in the higher dimensional Euclidean spaces with the Clifford algebra setting has attracted some researchers' attention (see e.g. \cite{Dang-Qian-Yang, Yang-Dang-Qian,Yang-Qian-Sommen,Dang-Qian-Chen}). In the study of (\ref{strong}) and (\ref{ex-strong}) we note that the derivative of $\varphi$ (or phase derivative) plays a role from the definitions of \lq\lq $\Cov_s$" and \lq\lq $\COV_s$". For a signal taking Clifford algebra values, one needs to define a proper phase derivative due to the non-commutativity of Clifford algebra. In \cite{Dang-Qian-Yang,Yang-Dang-Qian} several generalizations of phase derivative for Clifford algebra-valued signals are proposed. Correspondingly, analogous results of (\ref{ex-strong}) are also proved. Note that for periodic and spherical signals, there are parallel theories of uncertainty principle, where the periodic case was first studied by Breitenberger in \cite{Breitenberger} and the spherical case (i.e. functions defined on $\mathbb S^n\subset \mathbb R^{n+1}$) was generalized in e.g. \cite{Dai-Xu,Erb,Goh-Goodman,Narcowich,Rosler,Steinerberger} based on Breitenberger uncertainty principle. However, their results did not contain extra-strong uncertainty principle in those cases.  The extra-strong uncertainty principle for periodic signals was first studied in \cite{Dang}, and the spherical case was given only for $n=2$ due to the definition of phase derivative defined in \cite{Dang-Qian-Chen}.
	
	It would be interesting to consider whether there holds an extra-strong version of (a subclass of ) CKN inequalities for complex-valued and vector-valued functions. As mentioned previously, the key is to obtain the \lq\lq phase derivative" and \lq\lq COV" terms. To overcome these difficulties, in this paper we will propose a new notation, named \lq\lq generalized phase derivative",  that is a substitute of phase derivative in our settings, then we have a new formula for the \lq\lq COV" term to obtain an improved version of CKN inequalities for complex- and vector-valued functions. Moreover, using a similar technique we extend the extra-strong uncertainty principle given in \cite{Dang-Qian-Chen} to cases for $\mathbb {S}^n, n\geq 2$. 
	
	We first give the definition \lq\lq generalized phase derivative".
	\begin{defn}\label{generalized phase derivative} For $f=(f_1,f_2,...,f_m), f_j\in C^1(\mathbb R^n), j=1,2,\cdots,m,$ the generalized phase derivative is given by
		\begin{align}\label{phase-der-gen}
			|\Phi_f^\prime(x)|=\left(\sum_{k=1}^n\sum_{1\leq j<l\leq m}\left(\frac{f_l(x)}{|f(x)|^2}\partial_k f_j(x)-\frac{f_j(x)}{|f(x)|^2}\partial_k f_l(x)\right)^2\right)^\frac{1}{2},
		\end{align}
	where $\del = (\partial_1,...,\partial_n)^T = (\frac{\partial}{\partial x_1},...,\frac{\partial}{\partial x_n})^T.$ 
	
	For $f=(f_1,f_2,...,f_m), f_j\in C^1 (\mathbb S^n), j=1,2,\cdots,m,$ the generalized phase derivative is given by
		\begin{align}\label{phase-der-gen sphere}
			|\Phi_f^\prime(x)|=\left(\sum_{k=0}^n\sum_{1\leq j<l\leq m}\left(\frac{f_l(x)}{|f(x)|^2}\Gamma_k f_j(x)-\frac{f_j(x)}{|f(x)|^2}\Gamma_k f_l(x)\right)^2\right)^\frac{1}{2},
		\end{align}
		where  $\nabla_{\mathbb{S}^n}=(\Gamma_0,\cdots,\Gamma_n)^T$ is the gradient operator on sphere, and we introduce it in detail in \S 3 and for more properties of $\nabla_{\mathbb S^n}$ we refer to \cite{Goh-Goodman}.\end{defn}
		
		Note that generalized phase derivative can be considered as a generalization of phase derivative for vector-valued functions $f$. In fact, a complex-valued function $f$ could be written as $f(x)=u(x)+iv(x)=\rho(x)e^{i\varphi(x)}.$ Then, we have
		\begin{align*}
			\frac{\del f(x)}{f(x)}-\frac{\del \overline {f(x)}}{\overline {f(x)}} &= \frac{e^{i\varphi(x)}\del \rho(x)+i\rho(x)e^{i\varphi(x)}\del \varphi(x)}{\rho(x)e^{i\varphi(x)}}-\frac{e^{-i\varphi(x)}\del \rho(x)-i\rho(x)e^{-i\varphi(x)}\del \varphi(x)}{\rho(x)e^{-i\varphi(x)}}\\
			& = 2i\del \varphi(x).
		\end{align*}
		On the other hand, we can also write it as
		\begin{align*}
			\frac{\del f(x)}{f(x)}-\frac{\del \overline {f(x)}}{\overline {f(x)}} &=\frac{\overline {f(x)}\del f(x)}{|f(x)|^2}-\frac{f(x)\del \overline {f(x)}}{|f(x)|^2}\\
			& =\frac{(u(x)-iv(x))(\del u(x)+i\del v(x))-(u(x)+iv(x))(\del u(x)-i\del v(x))}{|f(x)|^2}\\
			& =2i\frac{u(x)\del v(x)-v(x)\del u(x)}{|f(x)|^2},
		\end{align*}
		which is used in the proof of Theorem \ref{improved-CKN-complex}.
		
		In this sense we can formally define $\Phi_f^\prime(x)$ as phase derivative for vector-valued functions $f$. For the unity of notation, we also use $\Phi_f^\prime(x)$ to denote phase derivative for complex-valued signals on $\mathbb R^n$ and $\mathbb S^n,$ that is, for $f=u+iv, u,v\in C^1(\mathbb R^n),$ we denote \[\Phi_f^\prime(x)=\frac{\del f(x)}{2if(x)}-\frac{\del \overline {f(x)}}{2i\overline {f(x)}}
	,\] and for $f=u+iv, u,v\in C^1(\mathbb S^n),$ 
		we denote \[\Phi_f^\prime(x)=\frac{\nabla_{\mathbb{S}^n} f(x)}{2if(x)}-\frac{\nabla_{\mathbb{S}^n} \overline {f(x)}}{2i\overline {f(x)}}.\]
		
		In the following we will abusively use the notation $\lq\lq${\rm COV}" to denote the absolute covariance in several cases although it has different formula in different case.

	Our main results are stated as follows. In \S 2, we obtain extra-strong version of CKN inequalities in several cases. 
	For complex-valued functions, we have
	\begin{theorem}\label{improved-CKN-complex}
		Suppose $f=u+iv, u,v\in C^\infty_0(\mathbb R^n),$   $0<q<2<p \text{ and }2<n<\frac{2(p-q)}{p-2}$.
		There holds
		\begin{align}
			\int_{\mathbb R^n} |\del f(x)|^2 dx \int_{\mathbb R^n}\frac{|f(x)|^{2p-2}}{|x|^{2q-2}} dx \geq \frac{(n-q)^2}{p^2}\left(\int_{\mathbb R^n} \frac{|f(x)|^p}{|x|^q}dx\right)^2 +\COV^2,
		\end{align}
		where
		\begin{align*}
			\COV=\int_{\mathbb R^n}\frac{|f(x)|^p}{|x|^{q-1}}\left|\Phi_f^\prime(x)\right| dx.
		\end{align*}
	\end{theorem}
	For vector-valued functions, we have
	\begin{theorem}\label{improved-CKN-vector}
		Suppose that $f=(f_1,f_2,...,f_m)$, and $f_j\in C^\infty_0(\mathbb R^n),1\leq j\leq m$ is real-valued and $0<q<2<p \text{ and }2<n<\frac{2(p-q)}{p-2}$.
		There holds
		\begin{align*}
			\int_{\mathbb R^n}|\del f(x)|^2 dx\int_{\mathbb R^n}\frac{|f(x)|^{2p-2}}{|x|^{2q-2}}dx\geq \frac{(n-q)^2}{p^2}\left(\int_{\mathbb R^n}\frac{|f(x)|^p}{|x|^q}dx\right)^2 +\COV^2,
		\end{align*}
		where
		\begin{align*}
			\COV = \int_{\mathbb R^n}\frac{|f(x)|^{p}}{|x|^{q-1}}|\Phi_f^\prime(x)|dx.
		\end{align*}
	\end{theorem}

From the proofs of Theorem \ref{improved-CKN-complex} and Theorem \ref{improved-CKN-vector}, we can immediately obtain the results for HPW uncertainty principle with $p=2, q=0.$
\begin{theorem}\label{extra-UP-complex}
	Suppose $f=u+iv, u,v\in C^\infty_0(\mathbb R^n).$
		There holds
		\begin{align*}
			\int_{\mathbb R^n} |\del f(x)|^2 dx \int_{\mathbb R^n}|x|^2|f(x)|^2 dx \geq \frac{n^2}{4}\left(\int_{\mathbb R^n} {|f(x)|^2}dx\right)^2 +\COV^2,
		\end{align*}
		where
		\begin{align*}
			\COV=\int_{\mathbb R^n}|f(x)|^2|x|\left|\Phi_f^\prime(x)\right| dx.
		\end{align*}\end{theorem}

	\begin{theorem}\label{extra-UP-vector}
		Suppose that $f=(f_1,f_2,...,f_m)$, and $f_j\in C^\infty_0(\mathbb R^n),1\leq j\leq m$ is real-valued.
		There holds
		\begin{align*}
			\int_{\mathbb R^n}|\del f(x)|^2 dx\int_{\mathbb R^n}{|x|^{2}}{|f(x)|^{2}}dx\geq \frac{n^2}{4}\left(\int_{\mathbb R^n}{|f(x)|^2}dx\right)^2 +\COV^2,
		\end{align*}
		where
		\begin{align*}
			\COV = \int_{\mathbb R^n}{|f(x)|^{2}}|x||\Phi_f^\prime(x)|dx.
		\end{align*}
	\end{theorem}

		As an application of Theorem \ref{improved-CKN-vector}, we have
	\begin{coro}\label{CKN_grad}
	Suppose that $f=\del u$, and $u\in C_0^\infty(\mathbb R^n)$ is real-valued, and $0<q<2<p$ and $2<n<\frac{2(p-q)}{p-2}.$ There holds
		\begin{align}\label{CKN-Rellich}
			\begin{split}
				&\int_{\mathbb R^n}|\Delta u(x)|^2dx \int_{\mathbb R^n}\frac{|\del u(x)|^{2p-2}}{|x|^{2q-2}}dx\\ &\geq \frac{(n-q)^2}{p^2}\left(\int_{\mathbb R^n}\frac{|\del u(x)|^p}{|x|^q}dx\right)^2
				+ \left(\int_{\mathbb R^n}\frac{|\del u(x)|^{p}}{|x|^{q-1}}|\Phi_{\del u}^\prime(x)| dx\right)^2.
			\end{split}
		\end{align}
		\end{coro}
	\begin{remark}		
		In \cite{Cazacu-Flynn-Lam} Cazacu, Flynn and Lam proved the second order uncertainty principle with a sharp constant, and more general results are given by Duong and Nguyen in \cite{Duong-Nguyen}. Our application as given in (\ref{CKN-Rellich}) also provides a second order CKN inequalities. Compared with the result given in \cite{Cazacu-Flynn-Lam}, the constant in (\ref{CKN-Rellich}) is not sharp but with a \lq\lq  \ \COV" term. On one hand, it would be interesting to find connection between Cazacu, Flynn and Lam's result and (\ref{CKN-Rellich}). On the other hand, one may expect a sharp second order uncertainty principle with \lq\lq\  \COV" term if combining these two technical results.
	\end{remark}
	For CKN inequalities with general parameters, we can give the following partial generalization.
	\begin{theorem}\label{improved-CKN-gen}
		Suppose $f=(f_1,f_2,...,f_m)$, and $f_j\in C^\infty_0(\mathbb R^n)$ is real-valued. Let $n\geq 2, r>p>2 $ and $\alpha,\beta, \gamma$ be fixed real numbers such that
		\begin{align*}
			\frac{1}{p}+\frac{\alpha}{n}>0,  \frac{p-1}{p(r-1)}+\frac{\beta}{n}>0,\frac{1}{r}+\frac{\gamma}{n}>0,
		\end{align*}
		and
		\begin{align*}
			\gamma=\frac{1}{r}(\alpha-1)+\frac{p-1}{pr}\beta.
		\end{align*}
		Then, there holds
		\begin{align*}
			&\left(\int_{\mathbb R^n}|x|^{\alpha p} |\del f(x)|^p dx\right)\left(\int_{\mathbb R^n}|x|^\beta |f(x)|^\frac{p(r-1)}{p-1}dx\right)^{p-1}\\
			\geq & \left(\frac{n+\gamma r}{r}\right)^p\left(\int_{\mathbb R^n} |x|^{\gamma r}|f(x)|^r dx\right)^p +  \left(\int_{\mathbb R^n}|x|^{r\gamma +1} |\Phi_f^\prime(x)| |f(x)|^{r}dx\right)^p.
		\end{align*}

	\end{theorem}

	In the last section, as an application of CKN inequalities with \lq\lq\COV" term, we obtain an extra-strong uncertainty principle for complex- and vector-valued signals on $n$-sphere $\mathbb S^n$, $n\geq 2$, which generalizes existing results for complex-valued signals on $n$-sphere $\mathbb S^n$, $n\geq 2$ \cite{Goh-Goodman}, for radial
functions \cite{Rosler}, for real-valued functions \cite{Narcowich} on $\mathbb S^2$ and for complex-valued functions \cite{Dang-Qian-Chen}. In the following, without loss of generality,  we assume that signals $f(x)$ on $\mathbb S^n$ are of unit energy, that is, $\int_{\mathbb S^n}|f(x)|^2d\sigma(x)=1.$ 
	To state our results on $\mathbb S^n,$ we introduce some notations.
	
%First we introduce the gradient operator proposed in \cite{Goh-Goodman}, which is given by %$\nabla_{\mathbb{S}^n}=(\Gamma_0,\cdots,\Gamma_n)^T.$ For more properties of $\nabla_{\mathbb S^n}$ %we refer to \cite{Goh-Goodman} (see also \S3).
\begin{defn}\label{defn mean on the sphere vector complex}Suppose that $f(x)\in C^1(\mathbb{S}^n)$ is complex-valued. Then the spherical mean,
or mean of the space vector variable $x,$ is defined to be
\begin{eqnarray}\label{mean of time vector complex} \tau_f \triangleq\int_{\mathbb{S}^n } x|f(x)|^2d\sigma(x),\end{eqnarray}
the variance of $x$ is defined to be
\begin{eqnarray}\label{variance of time vector complex}{\rm V}_{x,f}\triangleq\int_{\mathbb{S}^n}|x-\tau_f|^2
|f(x)|^2d\sigma(x)=1-|\tau_f|^2,\end{eqnarray}
the mean of
frequency is defined by the following two methods, respectively, given by
\begin{eqnarray}\label{mean of frequency vector complex}a(f)\triangleq\int_{\mathbb{S}^n}[\nabla_{\mathbb{S}^n}f(x)]
\overline{f(x)}d\sigma(x), \end{eqnarray} and\begin{eqnarray}\label{mean of frequency vector complex 2}a^*(f)\triangleq {\rm Im}\big\{\int_{\mathbb{S}^n}[\nabla_{\mathbb{S}^n}f(x)]
\overline{f(x)}d\sigma(x)\big\}, \end{eqnarray}
correspondingly, the variance of frequency has two formulations, defined by the following two formulas,
\begin{eqnarray}\label{variance of frequency vector complex}
{\rm V}_{\nabla_{\mathbb{S}^n},f}\triangleq
\int_{\mathbb{S}^n}|\nabla_{\mathbb{S}^n}f(x)-a(f) f(x)|^2d\sigma(x)
,\end{eqnarray}
and \begin{eqnarray}\label{variance of frequency vector complex 2}
{\rm V}_{\nabla_{\mathbb{S}^n},f}^*\triangleq
\int_{\mathbb{S}^n}|-i\nabla_{\mathbb{S}^n}f(x)-a^*(f) f(x)|^2d\sigma(x)
,\end{eqnarray}
and the covariance also has two definitions given by
\begin{eqnarray}\label{covariance para vector complex}
{\rm COV}\triangleq \int_{\mathbb{S}^n}|x-\tau_f|
|\Phi_f^\prime(x)||f(x)|^2d\sigma(x),
\end{eqnarray}
and \begin{eqnarray}\label{covariance para vector complex 2}
{\rm COV}_*\triangleq \int_{\mathbb{S}^n}|x-\tau_f|
|\Phi_f^\prime(x)-a^*(f)||f(x)|^2d\sigma(x),
\end{eqnarray}
provided that the above integrals are well-defined.
\end{defn}
\begin{remark}\label{two definitions of variances }
In Definition \ref{defn mean on the sphere vector complex},  the mean of frequency, variance of frequency and covariance all have two formulations, that is inspired by the definitions of those in the literature. For example, the mean of frequency for periodic signal $s(t)=\rho(t)e^{i\varphi(t)}\in L^2([0,2\pi))$ with the Fourier coefficient $c_k, k=0, \pm 1,\pm 2,\cdots$ is given by $\langle k\rangle \triangleq\sum_{k=-\infty}^{\infty}k|c_k|^2$  (see \cite{Dang,Prestin,Prestin2}), that can be represented in the time domain by Plancherel's Theorem as follow:
\begin{eqnarray*}\langle k\rangle =\sum_{k=-\infty}^{\infty}k|c_k|^2=-i\int_0^{2\pi}s'(t)\overline{s(t)}dt=\int_0^{2\pi}{\rm Im}\big[\frac{s'(t)}{s(t)}\big]|s(t)|^2dt=\int_0^{2\pi}\varphi'(t)|s(t)|^2dt.\end{eqnarray*} Since $s(t)$ is periodic, we know the real part of $\int_0^{2\pi}s'(t)\overline{s(t)}dt$ is zero, so $\langle k\rangle$ is also equal to $\int_0^{2\pi}s'(t)\overline{s(t)}dt.$ For signals on $\mathbb{S}^n,$ in the existing literature, mean of frequency also could be defined by two methods: one in time domain (see \cite{Dang-Qian-Yang,Goh-Goodman,Narcowich,Rosler}), and the other one in frequency domain (see \cite{Dang-Qian-Yang}). In the present paper, we use the gradient operator $\nabla_{\mathbb{S}^n}=(\Gamma_0,\cdots,\Gamma_n)^T$ to define mean of frequency. Since the real part of $\int_{\mathbb{S}^n}[\nabla_{\mathbb{S}^n}f(x)]
\overline{f(x)}d\sigma(x)$ is not zero, so we adopt two formulas (\ref{mean of frequency vector complex}) and (\ref{mean of frequency vector complex 2}) to define mean of frequency. The formula (\ref{mean of frequency vector complex 2}) is essentially corresponding to $\int_0^{2\pi}{\rm Im}\big[\frac{s'(t)}{s(t)}\big]|s(t)|^2dt$ for periodic signal. Accordingly, variance of frequency and covariance also have two formulas.
\end{remark}
	
\begin{theorem}\label{UP vector complex into two parts}
 Suppose that $f(x)$ is complex-valued with $f(x)$ and $|f(x)|\in C^1(\mathbb{S}^n).$ Then \begin{eqnarray}\label{UP vector complex into two parts formula introduction}
{\rm V}_{x,f}{\rm V}_{\nabla_{\mathbb{S}^n},f}&\geq&\frac{n^2}{4}|\tau_f|^4+{\rm COV}^2.
\end{eqnarray}
\end{theorem}

\begin{theorem}\label{UP 2 vector complex into two parts}
 Suppose that $f(x)$ is complex-valued with $f(x)$ and $|f(x)|\in C^1(\mathbb{S}^n).$ Then \begin{eqnarray}\label{UP 2 vector complex into two parts formula introduction}
{\rm V}_{x,f}{\rm V}_{\nabla_{\mathbb{S}^n},f}^*&\geq&\frac{n^2}{4}|\tau_f|^4+{\rm COV}_*^2.
\end{eqnarray}
\end{theorem}

\begin{coro}\label{UP-sphere-comp}
 Let $f(x)$ be complex-valued with $f(x)$ and $|f(x)|\in C^1(\mathbb{S}^n).$ Then 
 \begin{eqnarray}\label{UP vector complex into two parts formula introduction coro}
     {\rm V}_{x,f}\int_{\mathbb S^n}|\nabla_{\mathbb S^n}f(x)|^2d\sigma(x)\geq \frac{n^2}{4}|\tau_f|^2+\COV^2.
\end{eqnarray}
% or equivalently,
 %\begin{align}
 %\begin{split}
%&\left(1-\left|\int_{\mathbb S^n}x|f(x)|^2d\sigma(x)\right|^2\right)\int_{\mathbb %S^n}|\nabla_{\mathbb S^n}f(x)|^2d\sigma(x)\\
%\geq &\frac{n^2}{4}\left|\int_{\mathbb S^n}x|f(x)|^2d\sigma(x)\right|^2+{\rm %COV}^2.
%\end{split}
%\end{align}
\end{coro}
\begin{remark}\label{compare Goh and our }
In \cite{Goh-Goodman}, the authors considered the uncertainty principle for signals on $\mathbb S^n$, $n\geq 2.$ To compare their results and our results, assuming signals in \cite{Goh-Goodman} are of unit energy, we give the results in \cite{Goh-Goodman} in terms of our notations. The authors first obtained the uncertainty principle 
:\begin{eqnarray}\label{Goh UP}
{\rm V}_{x,f}{\rm V}_{\nabla_{\mathbb{S}^n},f}&\geq&\frac{n^2}{4}|\tau_f|^4,
\end{eqnarray}
then deduced the following corollary \begin{eqnarray}\label{Goh UP coro}
     {\rm V}_{x,f}\int_{\mathbb S^n}|\nabla_{\mathbb S^n}f(x)|^2d\sigma(x)\geq \frac{n^2}{4}|\tau_f|^2.
 \end{eqnarray}
In fact, it is important to have the corollary since it is the lower bound of (\ref{Goh UP coro}), which is the lower bound of the classical uncertainty principle on sphere. It is the corollary such that the uncertainty principle (\ref{Goh UP}) establishes the relationship with the classical uncertainty principle on sphere. 

 It is obvious that the lower bounds of (\ref{UP vector complex into two parts formula introduction}) and (\ref{UP vector complex into two parts formula introduction coro}) are greater that those of (\ref{Goh UP}) and (\ref{Goh UP coro}). So our results could be considered as a generalization of those in {Goh-Goodman}.
\end{remark}

For the vector-valued case, we abusively use the same notations as those in the complex-valued case.
\begin{defn}\label{defn mean on the sphere vector}Let $f(x)=(f_1, \cdots,f_m), f_j\in  C^1(\mathbb{S}^n), j=1,\cdots, m$ be real-valued with $\int_{\mathbb S^n}|f(x)|^2d\sigma(x)=1.$  Then the spherical mean,
or mean of the space vector variable $x,$ is defined to be
\begin{eqnarray}\label{mean of time vector} \tau_f \triangleq\int_{\mathbb{S}^n } x|f(x)|^2d\sigma(x),\end{eqnarray}
the variance of $x$ is defined to be
\begin{eqnarray}\label{variance of time vector}{\rm V}_{x,f}\triangleq\int_{\mathbb{S}^n}|x-\tau_f|^2
|f(x)|^2d\sigma(x)=1-|\tau_f|^2,\end{eqnarray}
the mean of
frequency is defined as
\begin{eqnarray}\label{mean of frequency vector}a(f)\triangleq\int_{\mathbb{S}^n}[\nabla_{\mathbb{S}^n}f(x)]
f^T(x)d\sigma(x), \end{eqnarray}
the variance of frequency is defined as
\begin{eqnarray}\label{variance of frequency vector}
{\rm V}_{\nabla_{\mathbb{S}^n},f}\triangleq
\int_{\mathbb{S}^n}|\nabla_{\mathbb{S}^n}f(x)-a(f) f(x)|^2d\sigma(x)
,\end{eqnarray}
and the covariance is defined by
\begin{eqnarray}\label{covariance para vector}
{\rm COV}^2\triangleq [\int_{\mathbb{S}^n}|x-\tau_f|
|\Phi'_f(x)||f(x)|^2d\sigma(x)]^2,
\end{eqnarray}
provided that the above integrals are well-defined, where we use the following notation for a $u\times v $ matrix $A$:  $|A|^2=\sum\limits_{j=1}^{u}\sum\limits_{k=1}^v a_{j,k}^2.$
\end{defn}	

	\begin{theorem}\label{UP vector into two parts}
 Let $f(x)=(f_1, \cdots,f_m), $ $f_j\in  C^1(\mathbb{S}^n), j=1,\cdots, m,$ Then \begin{eqnarray}\label{UP vector into two parts formula introduction}
{\rm V}_{x,f}{\rm V}_{\nabla_{\mathbb{S}^n},f}&\geq&\frac{n^2}{4}|\tau_f|^4+{\rm COV}^2,
\end{eqnarray}
which is equivalent to 
\begin{align}\label{UP-vector-2}
\begin{split}
{\rm V}_{x,f}\int_{\mathbb S^n}|\nabla_{\mathbb S^n}f(x)|^2d\sigma(x)
\geq  \frac{n^2}{4}|\tau_f|^2+{\rm COV}^2.
\end{split}
\end{align}
\end{theorem}
Unlike the complex-valued case, (\ref{UP vector into two parts formula introduction}) is not only sufficient but also necessary for (\ref{UP-vector-2}).
Similarly, Theorem \ref{UP vector into two parts} could be considered as a generalization of \cite[Corollary 5.1]{Goh-Goodman} and \cite[Theorem 4.4]{Dang-Qian-Chen} in the vector-valued case. In the forthcoming paper we will also generalize the results obtained in this paper to the setting of Clifford algebra.

The paper is organized as follows. In \S 2 some extra-strong CKN inequalities for complex- and vector-valued functions are proved. In \S 3 generalizations of uncertainty principle on $\mathbb S^n,n\geq 2,$ for the complex- and vector-valued functions are proved, respectively.

	\section{CKN inequalities: Proof of Theorems \ref{improved-CKN-complex}, \ref{improved-CKN-vector}, \ref{improved-CKN-gen} and Corollary \ref{CKN_grad}}
	In this section, we give the proof of sharper version of CKN inequalities in the complex-valued and vector-valued cases. 
	
	\subsection{The complex-valued case}
	\noindent\textbf{Proof of Theorem \ref{improved-CKN-complex}.} For any function $f=u+iv, u,v\in C^\infty_0(\mathbb R^n),$ we have
	\begin{align*}
		&\int_{\mathbb R^n} |\del f(x)|^2 dx \\
		=&\int_{\mathbb R^n} \langle \del f(x), \del f(x)\rangle dx\\
		=&\int_{\mathbb R^n} \langle \frac{\del f(x)}{f(x)},\frac{\del f(x)}{f(x)}\rangle |f(x)|^2 dx\\
		=&\int_{\mathbb R^n} \langle \frac{\overline {f(x)}\del f(x)}{|f(x)|^2},\frac{\overline {f(x)}\del f(x)}{|f(x)|^2}\rangle |f(x)|^2 dx\\
		=&\int_{\mathbb R^n}\langle \frac{(u(x)-iv(x))(\del u(x)+i\del v(x))}{|f(x)|^2},\frac{(u(x)-iv(x))(\del u(x)+i\del v(x))}{|f(x)|^2}\rangle |f(x)|^2 dx\\
		%&=\int_{\mathbb R^n}\langle \frac{(u(x)\del u(x)+v(x)\del v(x)+i(u(x)\del v(x)-v(x)\del u(x)))}{|f(x)|^2},\frac{(u(x)\del u(x)+v(x)\del v(x)+i(u(x)\del v(x)-v(x)\del u(x)))}{|f(x)|^2}\rangle |f(x)|^2 dx\\
		=&\int_{\mathbb R^n}\frac{|u(x)\del u(x)+v(x) \del v(x)|^2+|u(x)\del v(x)-v(x)\del u(x)|^2}{|f(x)|^4}|f(x)|^2dx\\
		=&\int_{\mathbb R^n}\frac{|u(x)\del u(x)+v(x)\del v(x)|^2}{|f(x)|^2}dx+\int_{\mathbb R^n}\frac{|u(x)\del v(x)-v(x)\del u(x)|^2}{|f(x)|^2}dx,
	\end{align*}
	where $\langle \cdot,\cdot\rangle$ denotes the complex Euclidean inner product of $\mathbb C^n.$
	
	\noindent Note that
	\begin{align*}
		\del |f(x)|^2=\del (u^2(x)+v^2(x))=2u(x)\del u(x)+2v(x)\del v(x),
	\end{align*}
	and we also have
	\begin{align*}
		\del |f(x)|^p=\frac{p}{2}|f(x)|^{p-2}\del |f(x)|^2.
	\end{align*}
	
	In the following we adopt the proof given in \cite[page 878]{Xia} and \cite[page 13]{Kri}. Then for any $x\in \mathbb R^n$ one has that
	\begin{align*}
		\Delta |x|^2 &=2n,\\
		|\del |x||&=1,\quad a.e.,\\
		2|x|\del |x|&=\del |x|^2.
	\end{align*}
	Using integration by parts and Cauchy-Schwartz's inequality, we have
	\begin{align*}
		&\int_{\mathbb R^n}\frac{|f(x)|^p}{|x|^q}dx\\
		=&\frac{1}{2n}\int_{\mathbb R^n}\frac{|f(x)|^p}{|x|^q}\Delta |x|^2 dx\\
		=&-\frac{1}{2n}\int_{\mathbb R^n}\langle \del |x|^2,\del\frac{|f(x)|^p}{|x|^q}\rangle dx\\
		=&-\frac{1}{2n}\int_{\mathbb R^n}\langle \del |x|^2,\del |f(x)|^p\rangle \frac{1}{|x|^q} dx -\frac{1}{2n}\int_{\mathbb R^n}\langle \del |x|^2,\del \left(\frac{1}{|x|^q}\right)\rangle |f(x)|^p dx\\
		=&-\frac{p}{2n}\int_{\mathbb R^n}\langle \del |x|,\frac{\del |f(x)|^2}{|f(x)|}\rangle \frac{|f(x)|^{p-1}}{|x|^{q-1}} dx +\frac{q}{n}\int_{\mathbb R^n}\langle \del |x|,\del |x|\rangle \frac{|f(x)|^p}{|x|^{q}} dx.
	\end{align*}
	Consequently, we have
	\begin{align*}
		\frac{n-q}{p}\int_{\mathbb R^n}\frac{|f(x)|^p}{|x|^q}dx \leq & \int_{\mathbb R^n}\left|\langle \del |x|, \frac{\del|f(x)|^2}{2|f(x)|} \rangle \right| \frac{|f|^{p-1}}{|x|^{q-1}} dx\\
		\leq &\int_{\mathbb R^n}| \del |x| \left| \frac{\gd|f(x)|^2}{2|f(x)|}\right| \frac{|f(x)|^{p-1}}{|x|^{q-1}} dx\\
		=&\int_{\mathbb R^n}\frac{|\del |f(x)|^2|}{2|f(x)|}\frac{|f(x)|^{p-1}}{|x|^{q-1}}dx\\
		\leq & \left(\int_{\mathbb R^n} \frac{|\del |f(x)|^2|^2}{4|f(x)|^2}dx \right)^\frac{1}{2}\left( \int_{\mathbb R^n}\frac{|f(x)|^{2p-2}}{|x|^{2q-2}}dx \right)^\frac{1}{2}\\
		=&\left(\int_{\mathbb R^n} \frac{|u(x)\del u(x)+v(x)\del v(x)|^2}{|f(x)|^2}dx \right)^\frac{1}{2}\left( \int_{\mathbb R^n}\frac{|f(x)|^{2p-2}}{|x|^{2q-2}}dx \right)^\frac{1}{2}.
	\end{align*}
	By Cauchy-Schwartz's inequality we have
	\begin{align*}
		&\int_{\mathbb R^n}\frac{|u(x)\del v(x)-v(x)\del u(x)|^2}{|f(x)|^2}dx \int_{\mathbb R^n}\frac{|f(x)|^{2p-2}}{|x|^{2q-2}}dx\\
		\geq & \left(\int_{\mathbb R^n}\frac{|u(x)\del v(x)-v(x)\del u(x)|}{|f(x)|}\frac{|f(x)|^{p-1}}{|x|^{q-1}}dx\right)^2\\
		= & \left(\int_{\mathbb R^n}\frac{|\overline {f(x)}\del f(x)-f(x)\del \overline {f(x)}|}{2|f(x)|^2}\frac{|f(x)|^{p}}{|x|^{q-1}}dx \right)^2\\
		= &\COV^2.
	\end{align*}
	Thus we obtain that
	\begin{align*}
		\int_{\mathbb R^n} |\del f(x)|^2 dx \int_{\mathbb R^n}\frac{|f(x)|^{2p-2}}{|x|^{2q-2}} dx \geq \frac{(n-q)^2}{p^2}\left(\int_{\mathbb R^n} \frac{|f(x)|^p}{|x|^q}dx\right)^2 +\COV^2.
	\end{align*}
	\endproof
\subsection{The vector-valued case}	
	\textbf{Proof of Theorem \ref{improved-CKN-vector}.}
	Note that $|\del f(x)|^2=\sum_{k=1}^n\sum_{j=1}^m |\partial_kf_j(x)|^2.$ The main trick is the following equality, i.e.,
	\begin{align*}
		&|\del f(x)|^2- |\del |f(x)||^2\\
		&=\sum_{k=1}^n\sum_{j=1}^m |\partial_kf_j(x)|^2- \frac{|\del |f(x)|^2|^2}{4|f(x)|^2}\\
		&=\sum_{k=1}^n\sum_{j=1}^m |\partial_kf_j(x)|^2- \sum_{k=1}^n |\sum_{j=1}^m\frac{f_j(x)}{|f(x)|}\partial_k f_j(x)|^2\\
		&=\sum_{k=1}^n\left(\sum_{j=1}^m |\partial_kf_j(x)|^2- \sum_{j=1}^m \frac{(f_j(x))^2}{|f(x)|^2}|\partial_kf_j(x)|^2-2\sum_{1\leq j<l\leq m}\frac{f_j(x)f_l(x)}{|f(x)|^2}\partial_k f_j(x)\partial_k f_l(x)\right)\\
		&=\sum_{k=1}^n\left(\sum_{j=1}^m |\partial_kf_j(x)|^2\left(1-\frac{(f_j(x))^2}{|f(x)|^2}\right)-2\sum_{1\leq j<l\leq m}\frac{f_j(x)f_l(x)}{|f(x)|^2}\partial_k f_j(x)\partial_k f_l(x)\right)\\
		&=\sum_{k=1}^n\sum_{1\leq j<l\leq m}\left(\frac{f_l(x)}{|f(x)|}\partial_k f_j(x)-\frac{f_j(x)}{|f(x)|}\partial_k f_l(x)\right)^2\\
		&= |f(x)|^2|\Phi_f^\prime(x)|^2,
	\end{align*}
	where $|\Phi_f^\prime(x)|$ is defined as (\ref{phase-der-gen}).
	Then we can obtain the desired result if we can show that
	\begin{align*}
		\int_{\mathbb R^n}|\del |f(x)||^2 dx\int_{\mathbb R^n}\frac{|f(x)|^{2p-2}}{|x|^{2q-2}}dx\geq \frac{(n-q)^2}{p^2}\left(\int_{\mathbb R^n}\frac{|f(x)|^p}{|x|^q}dx\right)^2,
	\end{align*}
	which immediately follows from the proof of Theorem \ref{improved-CKN-complex}. By Cauchy-Schwarz's inequality we have that
	\begin{align*}
		& \int_{\mathbb R^n} (|\del f(x)|^2-|\del |f(x)||^2) dx \int_{\mathbb R^n}\frac{|f(x)|^{2p-2}}{|x|^{2q-2}}dx \\
		= &\int_{\mathbb R^n}|f(x)|^2|\Phi_f^\prime(x)|^2 dx \int_{\mathbb R^n}\frac{|f(x)|^{2p-2}}{|x|^{2q-2}}dx \\
		\geq &\left(\int_{\mathbf R^{n}}|f(x)||\Phi_f^\prime(x)|\frac{|f(x)|^{p-1}}{|x|^{q-1}} dx\right)^2\\
		= & \left(\int_{\mathbf R^{n}}\frac{|f(x)|^{p}}{|x|^{q-1}}|\Phi_f^\prime(x)| dx\right)^2\\
		= & \COV^2.
	\end{align*}
	
	Therefore
	\begin{align*}
		\int_{\mathbb R^n} |\del f(x)|^2 dx \int_{\mathbb R^n}\frac{|f(x)|^{2p-2}}{|x|^{2q-2}}dx  \geq \frac{(n-q)^2}{p^2}\left(\int_{\mathbb R^n}\frac{|f(x)|^p}{|x|^q}dx\right)^2 +\COV^2.
	\end{align*}
	\endproof

\noindent\textbf{Proof of Corollary \ref{CKN_grad}.} First, for $f=(f_1,...,f_n),$ one has
\begin{align*}
    &|\sum_{j=1}^n\partial_jf_j(x)|^2 + \sum_{1\leq j<k\leq n}|\partial_jf_k(x)-\partial_kf_j(x)|^2\\
    =&\sum_{j=1}^n|\partial_jf_j(x)|^2+2\sum_{ 1\leq j<k\leq n}\partial_jf_j(x)\partial_kf_k(x)\\
	&+\sum_{ 1\leq j<k\leq n}(|\partial_jf_k(x)|^2+|\partial_kf_j(x)|^2-2\partial_jf_k(x)\partial_kf_j(x))\\
	=&\sum_{k=1}^n\sum_{j=1}^n|\partial_kf_j(x)|^2+2\sum_{ 1\leq j<k\leq n}\partial_jf_j(x)\partial_kf_k(x)-2\sum_{ 1\leq j<k\leq n}\partial_jf_k(x)\partial_kf_j(x).
\end{align*}

Now we let $f=\del u$. Then the above equality gives
\begin{align*}
    |\laplace u(x)|^2 = & |\laplace u(x)|^2 + \sum_{1\leq j<k\leq n}|\partial_j\partial_ku(x)-\partial_k\partial_ju(x)|^2\\
    = & \sum_{j=1}^n\sum_{k=1}^n|\partial_j\partial_ku(x)|^2 +2\sum_{1\leq j<k\leq n}\partial_j\partial_ju(x)\partial_k\partial_ku(x)\\
    &-2\sum_{1\leq j<k\leq n}\partial_j\partial_ku(x)\partial_k\partial_ju(x)\\
    = & |\del\del u(x)|^2+2\sum_{1\leq j<k\leq n}\partial_j\partial_ju(x)\partial_k\partial_ku(x)-2\sum_{1\leq j<k\leq n}\partial_j\partial_ku(x)\partial_k\partial_ju(x).
\end{align*}
Integrating to both sides, we have
\begin{align*}
    \int_{\mathbb R^n}|\laplace u(x)|^2 dx = & \int_{\mathbb R^n}|\del\del u(x)|^2dx +2\sum_{1\leq j<k\leq n}\int_{\mathbb R^n}\partial_j\partial_ju(x)\partial_k\partial_ku(x)dx\\
    &-2\sum_{1\leq j<k\leq n}\int_{\mathbb R^n}\partial_j\partial_ku(x)\partial_k\partial_ju(x)dx.
\end{align*}
Note that using integration by parts, we have
\begin{align*}
    \int_{\mathbb R^n}\partial_j\partial_ju(x)\partial_k\partial_ku(x)dx = -\int_{\mathbb R^n}\partial_ju(x)\partial_j\partial_k\partial_ku(x)dx= \int_{\mathbb R^n}\partial_k\partial_ju(x)\partial_j\partial_ku(x)dx,
\end{align*}
which implies that
\begin{align*}
    \int_{\mathbb R^n}|\del \del u(x)|^2dx = \int_{\mathbb R^n}|\laplace u(x)|^2 dx.
 \end{align*}
Then by Theorem \ref{improved-CKN-vector}, there holds
		\begin{align*}
				&\int_{\mathbb R^n}|\laplace u(x)|^2dx \int_{\mathbb R^n}\frac{|\del u(x)|^{2p-2}}{|x|^{2q-2}}dx\\
				%=&\int_{\mathbb R^n}|\del\del u(x)|^2dx \int_{\mathbb R^n}\frac{|\del u(x)|^{2p-2}}{|x|^{2q-2}}dx\\
				\geq &\frac{(n-q)^2}{p^2}\left(\int_{\mathbb R^n}\frac{|\del u(x)|^p}{|x|^q}dx\right)^2
				+ \left(\int_{\mathbb R^n}\frac{|\del u(x)|^{p}}{|x|^{q-1}}|\Phi_{\del u}^\prime(x)| dx\right)^2.
		\end{align*}
		\endproof

	\noindent\textbf{Proof of Theorem \ref{improved-CKN-gen}.}
	As shown in the proof of Theorem \ref{improved-CKN-vector}, we write
	\begin{align*}
		\int_{\mathbb R^n} |x|^{\alpha p}|\del f(x)|^p dx
		=\int_{\mathbb R^n}|x|^{\alpha p}\left(|\del|f(x)||^2 + |f(x)|^2|\Phi_f^\prime(x)|^2\right)^\frac{p}{2}dx.
	\end{align*}
	We also have
	\begin{align*}
		&\int_{\mathbb R^n}|x|^{\gamma r}|f(x)|^r dx\\
		&=\frac{r}{n}\int_{\mathbb R^n}\langle \del |x|,\frac{\del |f(x)|^2}{2|f(x)|}\rangle |x|^{\gamma r+1}|f(x)|^{r-1} dx -\frac{\gamma r}{n}\int_{\mathbb R^n}|x|^{\gamma r}|f(x)|^r dx,
	\end{align*}
	which implies that
	\begin{align*}
		\frac{n+\gamma r}{r}\int_{\mathbb R^n}|x|^{\gamma r}|f(x)|^r dx
		\leq \left(\int_{\mathbb R^n}|x|^{\alpha p}\frac{|\del |f(x)|^2|^p}{2^p|f(x)|^p}dx\right)^\frac{1}{p} \left(\int_{\mathbb R^n}|x|^\beta |f(x)|^\frac{p(r-1)}{p-1} dx\right)^\frac{p-1}{p}.
	\end{align*}

	To complete the proof, we use the following inequality
	\begin{align*}
		(a+b)^l\geq a^l+b^l, \quad \text{for $a\geq 0,b\geq 0$ and $l\geq 1$}.
	\end{align*}
	
	Then, for $p>2,$ we have
	\begin{align*}
		\int_{\mathbb R^n}|x|^{\alpha p}\left\{\left(|\del f(x)|^2\right)^\frac{p}{2}-\left(|\del |f(x)||^2\right)^\frac{p}{2}\right\}dx
		\geq \int_{\mathbb R^n}|x|^{\alpha p}|f(x)|^p|\Phi_f^\prime(x)|^pdx.
	\end{align*}
	Consequently, by H\"older's inequality we have
	\begin{align*}
		&\int_{\mathbb R^n}|x|^{\alpha p}\left\{|\del f(x)|^p-\left(|\del |f(x)||^2\right)^\frac{p}{2}\right\}dx \left(\int_{\mathbb R^n}|x|^\beta |f(x)|^\frac{p(r-1)}{p-1} dx\right)^{p-1}\\
		\geq & \int_{\mathbb R^n}|x|^{\alpha p}|f(x)|^p|\Phi_f^\prime(x)|^pdx \left(\int_{\mathbb R^n}|x|^\beta |f(x)|^\frac{p(r-1)}{p-1} dx\right)^{p-1}\\
		\geq & \left(\int_{\mathbb R^n}|x|^\alpha |f(x)||\Phi_f^\prime(x)||x|^\frac{\beta(p-1)}{p} |f(x)|^{r-1}dx\right)^p\\
		=& \left(\int_{\mathbb R^n}|x|^{\alpha +\frac{\beta(p-1)}{p}} |\Phi_f^\prime(x)| |f(x)|^{r}dx\right)^p\\
		=& \left(\int_{\mathbb R^n}|x|^{r\gamma +1} |\Phi_f^\prime(x)| |f(x)|^{r}dx\right)^p.
	\end{align*}
	Therefore,
	\begin{align*}
		&\left(\int_{\mathbb R^n}|x|^{\alpha p} |\del f(x)|^p dx\right)\left(\int_{\mathbb R^n}|x|^\beta |f(x)|^\frac{p(r-1)}{p-1}dx\right)^{p-1}\\
		\geq & \left(\frac{n+\gamma r}{r}\right)^p\left(\int_{\mathbb R^n} |x|^{\gamma r}|f(x)|^r dx\right)^p + \left(\int_{\mathbb R^n}|x|^{r\gamma +1} |\Phi_f^\prime(x)| |f(x)|^{r}dx\right)^p.
	\end{align*}

	\endproof

	\section{Uncertainty principle for sphere $\mathbb{S}^n$: Proof of Theorems \ref{UP vector complex into two parts}, \ref{UP 2 vector complex into two parts}, \ref{UP vector into two parts} and Corollary \ref{UP-sphere-comp} }
	In this section we use the generalized phase derivative \lq\lq $|\Phi_f^\prime(x)|$" introduced in Definition \ref{generalized phase derivative} to study uncertainty principle on $\mathbb S^n,n\geq 2.$ This gives a generalization of results proved in \cite{Dang-Qian-Chen} and \cite{Goh-Goodman}.
	
\subsection{The complex-valued case on $\mathbb S^n$}
Now we begin to consider uncertainty principle on  $\mathbb{S}^n=\{x=(x_0,\cdots,x_n)^T\in \mathbb{R}^{n+1}:|x|=1\}$ with normalized surface measure $\sigma.$ Let $L^2(\mathbb{S}^n)$ be the space of complex-valued square-integral functions on $\mathbb{S}^n$ with the inner product $\langle f, g\rangle:=\int_{\mathbb{S}^n}f\overline{g}d\sigma.$ In the paper, we adopt an operator $\nabla_{\mathbb{S}^n}=(\Gamma_0,\cdots,\Gamma_n)^T$ given in \cite{Goh-Goodman}. For the self-contained purpose, we give the definition of the operator in the following. Let $f$ be a complex-valued $C^1$ function on $\mathbb{S}^n$ and $x\in \mathbb{S}^n.$ Then the component of $\nabla_{\mathbb{S}^n}f(x)$ normal to the sphere at $x$ is zero, while any component of $\nabla_{\mathbb{S}^n}f(x)$ tangential to the sphere at $x$ is equal to the corresponding component of $\nabla f(x),$ i.e., $x\cdot \nabla_{\mathbb{S}^n}f(x)=0$ and for $y\in \mathbb{R}^{n+1}\backslash \{0\}$ with $y\cdot x=0,y\cdot\nabla_{\mathbb{S}^n}f(x)=y\cdot \nabla f(x).$ Now define $F: \mathbb{R}^{n+1}\backslash \{0\}\to \mathbb{R}$ by $F(x)=f(\frac{x}{|x|}).$ Since the component of $\nabla F(x)$ normal to the sphere at $x$ is zero, we have $\nabla_{\mathbb{S}^n}f(x)=\nabla F(x).$ In particular,this gives for $j,k=0,\cdots, n,$
\begin{equation}\label{nabla}
\Gamma_j(x_k) =\begin{cases}1-x_j^2,&j=k,\\
-x_jx_k,&j\neq k.\end{cases}\end{equation}	
	
In \cite{Goh-Goodman} the authors also derive an \lq\lq integration by parts" formula for the operator $\Gamma_j$, and we formulate it in the following lemma.

\begin{lemma}\label{integration by parts}
 Let $f(\underline{x})$ and $g\in C^1(\mathbb{S}^n).$ Then we have
  \begin{eqnarray}\label{integration by parts formula}\int_{\mathbb{S}^n}(\Gamma_j f)gd\sigma(x)=n\int_{\mathbb{S}^n}x_j f(x)g(x)d\sigma(x)-\int_{\mathbb{S}^n}f(\Gamma_j g)d\sigma(x). \end{eqnarray}

\end{lemma}

\begin{lemma}\label{variance of frequency vector complex into two parts}
 Let $f(x) $ be complex-valued with  $f(x)$ and $|f(x)|\in C^1(\mathbb{S}^n).$  Then the variance of frequency for $f(x) $ can be divided into two parts as the following two formulas give. \begin{align}\label{variance of frequency vector complex into two parts formula}
{\rm V}_{\nabla_{\mathbb{S}^n},f}=\int_{\mathbb{S}^n}|\nabla_{\mathbb{S}^n}|f(x)|-a(f)|f(x)||^2d\sigma(x)
+\int_{\mathbb{S}^n}|\Phi_f^\prime(x)|^2|f(x)|^2d\sigma(x),
\end{align}
and
\begin{align}\label{variance of frequency vector complex 2 into two parts formula}
{\rm V}_{\nabla_{\mathbb{S}^n},f}^*=\int_{\mathbb{S}^n}|\nabla_{\mathbb{S}^n}|f(x)||^2d\sigma(x)
+\int_{\mathbb{S}^n}|\Phi_f^\prime(x)-a^*(f)|^2|f(x)|^2d\sigma(x).
\end{align}
% \begin{align}\label{variance of frequency vector complex into two parts formula 2}
% {\rm V}_{\nabla_{\mathbb{S}^n},f}=\int_{\mathbb{S}^n}|\nabla_{\mathbb{S}^n}|f(x)||^2d\sigma(x)
% +\int_{\mathbb{S}^n}|\Phi_f^\prime(x)-a(f)|^2|f(x)|^2d\sigma(x)-\frac{n^2}{2}|\tau_f|^2.
% \end{align}
 \end{lemma}
\noindent{\textbf{Proof.}} 
Since $f(x) \in C^1(\mathbb{S}^n),$ we have $f(x)$ and $\nabla_{\mathbb{S}^n}f(x)\in L^2(\mathbb{S}^n),$ thus $a(f)$ and ${\rm V}_{\nabla_{\mathbb{S}^n},f}$ are well-defined. 

We first prove the formula (\ref{variance of frequency vector complex into two parts formula}). Set $G= \nabla_{\mathbb{S}^n}-a(f),$ as the proof of Theorem \ref{improved-CKN-complex}, we have
\begin{eqnarray*}
{\rm V}_{\nabla_{\mathbb{S}^n},f}=\int_{\mathbb{S}^n}|\nabla_{\mathbb{S}^n}f(x)-a(f) f(x)|^2d\sigma(x)
% \\
% &=&\int_{\mathbb{S}^n}\langle\frac{(u-iv)(Gu+iGv)}{|f(x)|^2}, \frac{(u-iv)(Gu+iGv)}{|f(x)|^2}\rangle|f(x)|^2 d\sigma(x)\\
% &=&
% \int_{\mathbb{S}^n}\langle\frac{uGu+vGv+i(u G v-v Gu)}{|f(x)|^2}, \frac{uGu+vGv+i(u G v-v Gu)}{|f(x)|^2}\rangle|f(x)|^2 d\sigma(x)\\
% &=&\int_{\mathbb{S}^n}\langle\frac{uGu+vGv}{|f(x)|^2}, \frac{uGu+vGv}{|f(x)|^2}\rangle|f(x)|^2 d\sigma(x)\\
% &&-i\int_{\mathbb{S}^n}\langle\frac{uGu+vGv}{|f(x)|^2}, \frac{u G v-v Gu}{|f(x)|^2}\rangle|f(x)|^2 d\sigma(x)\\
% &&+i\int_{\mathbb{S}^n}\langle\frac{u G v-v Gu}{|f(x)|^2}, \frac{uGu+vGv}{|f(x)|^2}\rangle|f(x)|^2 d\sigma(x)\\
% &&+\int_{\mathbb{S}^n}\langle\frac{u G v-v Gu}{|f(x)|^2}, \frac{u G v-v Gu}{|f(x)|^2}\rangle|f(x)|^2 d\sigma(x)\\&=&\int_{\mathbb{S}^n}\langle\frac{uGu+vGv}{|f(x)|^2}, \frac{uGu+vGv}{|f(x)|^2}\rangle|f(x)|^2 d\sigma(x)\\
% &&+\int_{\mathbb{S}^n}\langle\frac{u G v-v G u}{|f(x)|^2}, \frac{u G v-v Gu}{|f(x)|^2}\rangle|f(x)|^2 d\sigma(x)
=
\int_{\mathbb{S}^n}[\frac{|uGu+vGv|^2}{|f(x)|^2}
+\frac{|u G v-v Gu|^2}{|f(x)|^2}]d\sigma(x),
\end{eqnarray*}
which gives
\begin{eqnarray*}
{\rm V}_{\nabla_{\mathbb{S}^n},f}&=&
\int_{\mathbb{S}^n}|\nabla_{\mathbb{S}^n}|f(x)|-a(f)|f(x)||^2d\sigma(x)\\&&
+\int_{\mathbb{S}^n}\frac{|\overline{f(x)}[\nabla_{\mathbb{S}^n}f(x)-a(f)f(x)]-f(x)
[\nabla_{\mathbb{S}^n}\overline{f(x)}-a(f)\overline{f(x)}]|^2}{4|f(x)|^2}]d\sigma(x)\\
&=&
\int_{\mathbb{S}^n}|\nabla_{\mathbb{S}^n}|f(x)|-a(f)|f(x)||^2d\sigma(x)
+\int_{\mathbb{S}^n}\frac{|\overline{f(x)}\nabla_{\mathbb{S}^n}f(x)-f(x)
\nabla_{\mathbb{S}^n}\overline{f(x)}|^2}{4|f(x)|^2}d\sigma(x).
\end{eqnarray*}
Now it is time to prove (\ref{variance of frequency vector complex 2 into two parts formula}).
\begin{eqnarray*}
{\rm V}_{\nabla_{\mathbb{S}^n},f}^*&=&\int_{\mathbb{S}^n}|-i\nabla_{\mathbb{S}^n}f(x)-a^*(f) f(x)|^2d\sigma(x)\\
% &=&\int_{\mathbb{S}^n}\langle -i\nabla_{\mathbb{S}^n}f(x)-a^*(f)f(x),-i\nabla_{\mathbb{S}^n}f(x)-a^*(f)f(x) \rangle d\sigma(x)\\
% &=&\int_{\mathbb{S}^n}|\nabla_{\mathbb{S}^n}f(x)|^2 d\sigma(x)+\int_{\mathbb{S}^n}|a^*(f)|^2|f(x)|^2 d\sigma(x)\\&&+\int_{\mathbb{S}^n}\langle -i\nabla_{\mathbb{S}^n}f(x),-a^*(f)f(x) \rangle d\sigma(x)+\int_{\mathbb{S}^n}\langle -a^*(f)f(x),-i\nabla_{\mathbb{S}^n}f(x) \rangle d\sigma(x)\\&=&\int_{\mathbb{S}^n}|\nabla_{\mathbb{S}^n}f(x)|^2 d\sigma(x)+|a^*(f)|^2\\&&+i\sum_{j=0}^n a_j^*(f)\int_{\mathbb{S}^n}
% \Gamma_j f(x)\overline{f(x)}  d\sigma(x)-i\sum_{j=0}^na_j^*(f)\int_{\mathbb{S}^n}f(x)\overline{\Gamma_j f(x)} d\sigma(x)\\&=&\int_{\mathbb{S}^n}|\nabla_{\mathbb{S}^n}f(x)|^2 d\sigma(x)+|a^*(f)|^2\\&&+i\sum_{j=0}^n a_j^*(f)[\int_{\mathbb{S}^n}
% \Gamma_j f(x)\overline{f(x)}  d\sigma(x)-\int_{\mathbb{S}^n}f(x)\overline{\Gamma_j f(x)} d\sigma(x)]\\&=&\int_{\mathbb{S}^n}|\nabla_{\mathbb{S}^n}f(x)|^2 d\sigma(x)+|a^*(f)|^2+i\sum_{j=0}^n a_j^*(f)2ia_j^*(f)\\&=&\int_{\mathbb{S}^n}|\nabla_{\mathbb{S}^n}f(x)|^2 d\sigma(x)+|a^*(f)|^2-2|a^*(f)|^2\\
&=&\int_{\mathbb{S}^n}|\nabla_{\mathbb{S}^n}f(x)|^2 d\sigma(x)-|a^*(f)|^2\\
%&=&\int_{\mathbb{S}^n}[\frac{|u\nabla_{\mathbb{S}^n}u+v\nabla_{\mathbb{S}^n}v|^2}{|f(x)|^2}
%+\frac{|u \nabla_{\mathbb{S}^n} v-v \nabla_{\mathbb{S}^n}u|^2}{|f(x)|^2}]d\sigma(x)-|a(f)|^2\\
&=&\int_{\mathbb{S}^n}|\nabla_{\mathbb{S}^n}|f(x)||^2 d\sigma(x)+\int_{\mathbb{S}^n}\frac{|\overline{f(x)}\nabla_{\mathbb{S}^n}f(x)-f(x)
\nabla_{\mathbb{S}^n}\overline{f(x)}|^2}{4|f(x)|^2}d\sigma(x)-|a^*(f)|^2\\
&=&\int_{\mathbb{S}^n}|\nabla_{\mathbb{S}^n}|f(x)||^2 d\sigma(x)+\int_{\mathbb{S}^n}|\frac{\overline{f(x)}\nabla_{\mathbb{S}^n}f(x)-f(x)
\nabla_{\mathbb{S}^n}\overline{f(x)}}{2i|f(x)|^2}-a^*(f)|^2|f(x)|^2d\sigma(x),
\end{eqnarray*}
where in the last equality we use the formula
\begin{eqnarray}\label{equality in proof} &&\int_{\mathbb{S}^n}|\frac{\overline{f(x)}\nabla_{\mathbb{S}^n}f(x)-f(x)
\nabla_{\mathbb{S}^n}\overline{f(x)}}{2i|f(x)|^2}-a^*(f)|^2|f(x)|^2d\sigma(x)\nonumber\\&=&\int_{\mathbb{S}^n}|\frac{\overline{f(x)}\nabla_{\mathbb{S}^n}f(x)-f(x)
\nabla_{\mathbb{S}^n}\overline{f(x)}}{2\textbf{}|f(x)|^2}|^2|f(x)|^2d\sigma(x)-|a^*(f)|^2.
\end{eqnarray}
Now we prove the formula (\ref{equality in proof}).
\begin{eqnarray*} 
&&\int_{\mathbb{S}^n}|\frac{\overline{f(x)}\nabla_{\mathbb{S}^n}f(x)-f(x)
\nabla_{\mathbb{S}^n}\overline{f(x)}}{2i|f(x)|^2}-a^*(f)|^2|f(x)|^2d\sigma(x)\\
&=&\int_{\mathbb{S}^n}\langle\frac{\overline{f(x)}\nabla_{\mathbb{S}^n}f(x)-f(x)
\nabla_{\mathbb{S}^n}\overline{f(x)}}{2i|f(x)|^2}-a^*(f),\frac{\overline{f(x)}\nabla_{\mathbb{S}^n}f(x)-f(x)
\nabla_{\mathbb{S}^n}\overline{f(x)}}{2i|f(x)|^2}-a^*(f)\rangle|f(x)|^2d\sigma(x)\\
&=&\int_{\mathbb{S}^n}|\frac{\overline{f(x)}\nabla_{\mathbb{S}^n}f(x)-f(x)
\nabla_{\mathbb{S}^n}\overline{f(x)}}{2|f(x)|^2}|^2|f(x)|^2d\sigma(x)+\int_{\mathbb{S}^n}|a^*(f)|^2|f(x)|^2d\sigma(x)\\&&-
\int_{\mathbb{S}^n}\langle\frac{\overline{f(x)}\nabla_{\mathbb{S}^n}f(x)-f(x)
\nabla_{\mathbb{S}^n}\overline{f(x)}}{2i|f(x)|^2},a^*(f)\rangle|f(x)|^2d\sigma(x)\\
&&-\int_{\mathbb{S}^n}\langle a^*(f),\frac{\overline{f(x)}\nabla_{\mathbb{S}^n}f(x)-f(x)
\nabla_{\mathbb{S}^n}\overline{f(x)}}{2i|f(x)|^2}\rangle|f(x)|^2d\sigma(x)\\
&=&\int_{\mathbb{S}^n}|\frac{\overline{f(x)}\nabla_{\mathbb{S}^n}f(x)-f(x)
\nabla_{\mathbb{S}^n}\overline{f(x)}}{2|f(x)|^2}|^2|f(x)|^2d\sigma(x)+|a^*(f)|^2-2|a^*(f)|^2\\
&=&\int_{\mathbb{S}^n}|\frac{\overline{f(x)}\nabla_{\mathbb{S}^n}f(x)-f(x)
\nabla_{\mathbb{S}^n}\overline{f(x)}}{2|f(x)|^2}|^2|f(x)|^2d\sigma(x)-|a^*(f)|^2.
\end{eqnarray*}

\endproof

%\begin{theorem}\label{UP vector complex into two parts}
% Let $f(x)$ be complex-valued with $f(x)$ and $xf(x)\in L^2(\mathbb{S}^n)$,  all the first order partial derivatives of $f(x)$ and $|f(x)|$ exist and continuous, and $\nabla_{\mathbb{S}^n}f(x)\in L^2(\mathbb{S}^n).$ Then \begin{eqnarray}\label{UP vector complex into two parts formula}
%{\rm V}_{x,f}{\rm V}_{\nabla_{\mathbb{S}^n},f}&\geq&\frac{n^2}{4}|\tau_f|^4+{\rm COV}^2.
%\end{eqnarray}
%\end{theorem}
\noindent{\textbf{Proof of Theorem \ref{UP vector complex into two parts}.}}
By (\ref{variance of frequency vector complex into two parts formula}), we first prove that \begin{eqnarray}\label{UP complex formula one in proof}{\rm V}_{x,f}
\int_{\mathbb{S}^n}|\nabla_{\mathbb{S}^n}|f(x)|-a(f)|f(x)||^2d\sigma(x)\geq \frac{n^2}{4}|\tau_f|^4.\end{eqnarray} By Cauchy-Schwarz's inequality, we have
\begin{eqnarray*}
&&{\rm V}_{x,f}
\int_{\mathbb{S}^n}|\nabla_{\mathbb{S}^n}|f(x)|-a(f)|f(x)||^2d\sigma(x)\\
&=&\int_{\mathbb{S}^n}|x-\tau_f|^2
|f(x)|^2d\sigma(x)
\int_{\mathbb{S}^n}|\nabla_{\mathbb{S}^n}|f(x)|-a(f)|f(x)||^2d\sigma(x)\\
&\geq&[\int_{\mathbb{S}^n}|x-\tau_f||f(x)|
|\nabla_{\mathbb{S}^n}|f(x)|-a(f)|f(x)||d\sigma(x)]^2\\
&\geq&[\int_{\mathbb{S}^n}
(\nabla_{\mathbb{S}^n}|f(x)|-a(f)|f(x)|)\cdot(x|f(x)|-\tau_f|f(x)|)d\sigma(x)]^2\\
&=&\{\int_{\mathbb{S}^n}
[(\nabla_{\mathbb{S}^n}|f(x)|)\cdot(x|f(x)|)
-(\nabla_{\mathbb{S}^n}|f(x)|)\cdot (\tau_f|f(x)|)\\&&
-(a(f)|f(x)|)\cdot(x|f(x)|)+(a(f)|f(x)|)\cdot(\tau_f|f(x)|)]d\sigma(x)\}^2\\
&=&\{\int_{\mathbb{S}^n}
[(\nabla_{\mathbb{S}^n}|f(x)|)\cdot(x|f(x)|)
-(\nabla_{\mathbb{S}^n}|f(x)|)\cdot (\tau_f|f(x)|)]d\sigma(x)\}^2\\
&=&\{\int_{\mathbb{S}^n}
\sum_{j=0}^n[(\Gamma_j|f(x)|)(x_j|f(x)|)
-(\Gamma_j|f(x)|) (\tau_{f,j}|f(x)|)]d\sigma(x)\}^2\\
&=&\{\sum_{j=0}^n\int_{\mathbb{S}^n}
[(\Gamma_j|f(x)|)(x_j|f(x)|)
-(\Gamma_j|f(x)|) (\tau_{f,j}|f(x)|)]d\sigma(x)\}^2\\&=&A,\end{eqnarray*}
where we denote $\tau_{f}=(\tau_{f,0},\tau_{f,1},\cdots, \tau_{f,n})^T.$

Now we further calculate $A.$
By Lemma \ref{integration by parts} and (\ref{nabla}), we have
\begin{eqnarray*}
&&\int_{\mathbb{S}^n}
 (\Gamma_j|f(x)|)(x_j|f(x)|)d\sigma(x)\\
 &=&n\int_{\mathbb{S}^n}
 x_j|f(x)|x_j|f(x)|d\sigma(x)-\int_{\mathbb{S}^n}
|f(x)|[\Gamma_jx_j|f(x)|]d\sigma(x)\\
 &=&n\int_{\mathbb{S}^n}
 x_j^2|f(x)|^2d\sigma(x)-\{\int_{\mathbb{S}^n}
|f(x)|[\Gamma_jx_j]|f(x)|d\sigma(x)+\int_{\mathbb{S}^n}
|f(x)|x_j[\Gamma_j|f(x)|]d\sigma(x)\}\\
 &=&n\int_{\mathbb{S}^n}
 x_j^2|f(x)|^2d\sigma(x)-\int_{\mathbb{S}^n}
|f(x)|[1-x_j^2]|f(x)|d\sigma(x)-\int_{\mathbb{S}^n}
|f(x)|x_j[\Gamma_j|f(x)|]d\sigma(x),
\end{eqnarray*}then
\begin{eqnarray}\label{part1}
\int_{\mathbb{S}^n}
 (\Gamma_j|f(x)|)(x_j|f(x)|)d\sigma(x)=\frac{n+1}{2}\int_{\mathbb{S}^n}
 x_j^2|f(x)|^2d\sigma(x)-\frac{1}{2}.
\end{eqnarray}
By Lemma \ref{integration by parts} again,
\begin{eqnarray*}
\int_{\mathbb{S}^n}
 (\Gamma_j|f(x)|) (\tau_{f,j}|f(x)|)d\sigma(x)=\tau_{f,j}\{n\int_{\mathbb{S}^n}
 x_j|f(x)||f(x)|d\sigma(x)-\int_{\mathbb{S}^n}
|f(x)|[\Gamma_j|f(x)|]d\sigma(x)\},
\end{eqnarray*}then we have
\begin{eqnarray}\label{part2}
\int_{\mathbb{S}^n}
 (\Gamma_j|f(x)|) (\tau_{f,j}|f(x)|)d\sigma(x)=\frac{n}{2}\tau_{f,j}^2.
\end{eqnarray}
Inserting (\ref{part1}) and (\ref{part2}), we obtain
\begin{eqnarray*}
&&{\rm V}_{x,f}
\int_{\mathbb{S}^n}|\nabla_{\mathbb{S}^n}|f(x)|-a(f)|f(x)||^2d\sigma(x)\\&\geq&\{\sum_{j=0}^n[\frac{n+1}{2}\int_{\mathbb{S}^n}
 x_j^2|f(x)|^2d\sigma(x)-\frac{1}{2}+\frac{n}{2}\tau_{f,j}^2]\}^2\\&=&\frac{n^2}{4}|\tau_f|^4.
\end{eqnarray*}
Now we prove \[{\rm V}_{x,f}
\int_{\mathbb{S}^n}|\Phi_f^\prime(x)|^2|f(x)|^2d\sigma(x)\geq {\rm COV}^2.\]

By using Cauchy-Schwarz's inequality, we immediately obtain
\begin{eqnarray*}
&&{\rm V}_{x,f}
\int_{\mathbb{S}^n}|\Phi_f^\prime(x)|^2|f(x)|^2d\sigma(x)\\
&=&\int_{\mathbb{S}^n}|x-\tau_f|^2
|f(x)|^2d\sigma(x)
\int_{\mathbb{S}^n}|\Phi_f^\prime(x)|^2|f(x)|^2d\sigma(x)\\
&\geq&[\int_{\mathbb{S}^n}|x-\tau_f|
|f(x)||\Phi_f^\prime(x)||f(x)|d\sigma(x)]^2\\
&=& \COV^2.
\end{eqnarray*}
\endproof

\noindent{\textbf{Proof of Theorem \ref{UP 2 vector complex into two parts}.}}
The proof is similar with that of Theorem \ref{UP vector complex into two parts}.
\endproof

\medskip

One can obtain Corollary \ref{UP-sphere-comp} by the proof of \cite{Goh-Goodman} immediately. For the self-contained purpose, we include a proof.

\noindent{\textbf{Proof of Corollary \ref{UP-sphere-comp}.}}	
First, by the last formula in the proof of Lemma \ref{variance of frequency vector complex into two parts}, we have 
\begin{align*}
    {\rm Re}\{\int_{\mathbb S^n}(\Gamma_jf(x))\overline{f(x)}d\sigma(x)\} =\frac{n}{2}\int_{\mathbb S^n}x_j|f(x)|^2 d\sigma(x)
\end{align*}
and consequently,
\begin{align}\label{af_tauf}
    \left|\int_{\mathbb S^n}(\nabla_{\mathbb S^n}f(x))\overline{f(x)}d\sigma(x)\right|^2 \geq \sum_{j=0}^n\frac{n^2}{4}\left|\int_{\mathbb S^n}x_j|f(x)|^2d\sigma(x)\right|^2=\frac{n^2}{4}\left|\int_{\mathbb S^n}x|f(x)|^2d\sigma(x)\right|^2.
\end{align}
By Theorem \ref{UP vector complex into two parts}, we have
\begin{align*}
    &\int_{\mathbb S^n}|\nabla_{\mathbb S^n}f(x)-a(f)f(x)|^2d\sigma(x)\left(1-\left|\int_{\mathbb S^n}x|f(x)|^2d\sigma(x)\right|^2\right)\\
    = & \left(\int_{\mathbb S^n}|\nabla_{\mathbb S^n}f(x)|^2d\sigma(x)-|a(f)|^2\right)\left(1-|\tau_f|^2\right)\\\geq &  \frac{n^2}{4}|\tau_f|^4 +\COV^2,
\end{align*}
which yields
\begin{align*}
&\int_{\mathbb S^n}|\nabla_{\mathbb S^n}f(x)|^2d\sigma(x)\left(1-|\tau_f|^2\right)\\
\geq & \frac{n^2}{4}|\tau_f|^4 +|a(f)|^2-|a(f)|^2|\tau_f|^2+\COV^2\\
    \geq & \frac{n^2}{4}|\tau_f|^4 +|a(f)|^2-|a(f)|^2|\tau_f|^2+\COV^2+\frac{n^2}{4}|\tau_f|^2 -\frac{n^2}{4}|\tau_f|^2\\
    \geq & \frac{n^2}{4}|\tau_f|^2\left(|\tau_f|^2-1\right)+\COV^2+\frac{n^2}{4}|\tau_f|^2 +|a(f)|^2\left(1-|\tau_f|^2\right)\\
    \geq & \left(1-|\tau_f|^2\right)\left(|a(f)|^2-\frac{n^2}{4}|\tau_f|^2\right) +\COV^2 +\frac{n^2}{4}|\tau_f|^2\\
    \geq & \COV^2 +\frac{n^2}{4}|\tau_f|^2,
\end{align*}
where we have used the fact that (\ref{af_tauf}) and $|\tau_f|^2<1.$
\endproof

\subsection{The vector-valued case on sphere }

Now we begin to consider uncertainty principle for vector-valued $f=(f_1, \cdots,f_m)$ on $\mathbb{S}^n$ with normalized surface measure $\sigma.$ For convenience,  we abusively use the same notations as those in the complex-valued case.

\begin{lemma}\label{variance of frequency vector into two parts}
 Let $f(x)=(f_1, \cdots,f_m), f_j\in C^1(\mathbb{S}^n), j=1,\cdots, m$ be real-valued with $\int_{\mathbb S^n}|f(x)|^2d\sigma(x)=1.$ Then \begin{eqnarray}\label{variance of frequency vector into two parts formula}
{\rm V}_{\nabla_{\mathbb{S}^n},f}=\int_{\mathbb{S}^n}|\nabla_{\mathbb{S}^n}|f(x)|-a(f)|f(x)||^2d\sigma(x)
+\int_{\mathbb{S}^n}|\Phi'_f(x)|^2|f(x)|^2d\sigma(x),
\end{eqnarray}
where we recall that
\begin{align*}
 |\Phi'_f(x)|^2=  \sum\limits_{i=0}^{n}\sum\limits_{1\leq j<k \leq m}\big\{\frac{f_k(x)\Gamma_i f_j(x)}{|f(x)|^2}-\frac{f_j(x)\Gamma_i f_k(x)}{|f(x)|^2}\big\}^2.
\end{align*}
\end{lemma}
\noindent{\textbf{Proof.}}
Since $f_j\in C^1(\mathbb{S}^n), j=1,\cdots, m,$ we have $f_j$ and $\Gamma_if_j\in L^2(\mathbb{S}^n)$ for $i=0,1,\cdots,n, j=1,\cdots, m,$   then $a(f)$ and ${\rm V}_{\nabla_{\mathbb{S}^n},f}$ are well-defined. $a(f)$ is $(n+1)$-vector-valued and denoted as \[a(f)=(a_0(f),a_1(f),\cdots, a_n(f))^T.\]
By direct computation we have \begin{eqnarray*}|\nabla_{\mathbb{S}^n}f(x)-a(f) f(x)|^2=\sum\limits_{i=0}^{n}\sum\limits_{j=1}^m[\Gamma_if_j-a_i(f)f_j]^2
\end{eqnarray*}
and \begin{eqnarray*}&&|\nabla_{\mathbb{S}^n}|f(x)|-a(f) |f(x)||^2\\&=&\sum\limits_{i=0}^{n}[\Gamma_i|f(x)|-a_i(f)|f(x)|]^2\\
&=&\sum\limits_{i=0}^{n}\frac{[\Gamma_i|f(x)|^2-2a_i(f)|f(x)|^2]^2}{4|f(x)|^2}\\
&=&\sum\limits_{i=0}^{n}\frac{\big\{\Gamma_i[\sum\limits_{j=1}^m f_j^2(x)]-2a_i(f)[\sum\limits_{j=1}^m f_j^2(x)]\big\}^2}{4|f(x)|^2}\\
&=&\sum\limits_{i=0}^{n}\frac{\big\{\sum\limits_{j=1}^m[ f_j(x)\Gamma_i f_j(x)-a_i(f) f_j^2(x)]\big\}^2}{|f(x)|^2}\\&=&\sum\limits_{i=0}^{n}\big\{\sum\limits_{j=1}^m \frac{[f_j(x)\Gamma_i f_j(x)-a_i(f) f_j^2(x)]^2}{|f(x)|^2}\\&&+2\sum\limits_{1\leq j<k \leq m}\frac{[f_j(x)\Gamma_i f_j(x)-a_i(f) f_j^2(x)][f_k(x)\Gamma_i f_k(x)-a_i(f) f_k^2(x)]}{|f(x)|^2}\big\}.\end{eqnarray*}

Furthermore, we can have
\begin{eqnarray*}&&|\nabla_{\mathbb{S}^n}f(x)-a(f) f(x)|^2-|\nabla_{\mathbb{S}^n}|f(x)|-a(f) |f(x)||^2\\&=&\sum\limits_{i=0}^{n}\sum\limits_{j=1}^m[\Gamma_if_j-a_i(f)f_j]^2-\sum\limits_{i=0}^{n+1}\sum\limits_{j=1}^m \frac{[f_j(x)\Gamma_i f_j(x)-a_i(f) f_j^2(x)]^2}{|f(x)|^2}\\&&
-2\sum\limits_{i=0}^{n+1}\sum\limits_{1\leq j<k \leq m}\frac{[f_j(x)\Gamma_i f_j(x)-a_i(f) f_j^2(x)][f_k(x)\Gamma_i f_k(x)-a_i(f) f_k^2(x)]}{|f(x)|^2}\\
&=&\sum\limits_{i=0}^{n}\sum\limits_{j=1}^m[\Gamma_if_j-a_i(f)f_j]^2(1-\frac{f_j^2(x)}{|f(x)|^2})\\
&&
-2\sum\limits_{i=0}^{n}\sum\limits_{1\leq j<k \leq m}\frac{f_j(x)f_k(x)[\Gamma_i f_j(x)-a_i(f) f_j(x)][\Gamma_i f_k(x)-a_i(f) f_k(x)]}{|f(x)|^2}\\
&=&\sum\limits_{i=0}^{n}\sum\limits_{1\leq j<k \leq m}\big\{\frac{f_k(x)[\Gamma_i f_j(x)-a_i(f) f_j(x)]}{|f(x)|}-\frac{f_j(x)[\Gamma_i f_k(x)-a_i(f) f_k(x)]}{|f(x)|}\big\}^2\\
&=&\sum\limits_{i=0}^{n}\sum\limits_{1\leq j<k \leq m}\big[\frac{f_k(x)}{|f(x)|}\Gamma_i f_j(x) -\frac{f_j(x)}{|f(x)|}\Gamma_i f_k(x) \big]^2\\
&=&|\Phi'_f(x)|^2|f(x)|^2.
\end{eqnarray*}
\endproof

%\begin{theorem}\label{UP vector into two parts}
% Let $f(x)=(f_1, \cdots,f_m), $ $x_if(x)$ and $f_j\in L^2(\mathbb{S}^n),$ all the first order partial derivatives of $f_j(x)$ and $|f_j(x)|$ exist and are continuous, and $\Gamma_if_j\in L^2(\mathbb{S}^n), i=0, 1,\cdots, n, j=1,\cdots, m,$ Then \begin{eqnarray}\label{UP vector into two parts formula}
%{\rm V}_{x,f}{\rm V}_{\nabla_{\mathbb{S}^n},f}&\geq&\frac{n^2}{4}|\tau_f|^4+{\rm COV}^2.
%\end{eqnarray}
%\end{theorem} \ref{UP vector complex into two parts} Corollary \ref{UP-sphere-comp} \
\noindent{\textbf{Proof of Theorem \ref{UP vector into two parts}.}}
By (\ref{variance of frequency vector into two parts formula}), we first prove that \[{\rm V}_{x,f}
\int_{\mathbb{S}^n}|\nabla_{\mathbb{S}^n}|f(x)|-a(f)|f(x)||^2d\sigma(x)\geq \frac{n^2}{4}|\tau_f|^4,\] which can be proved by the same way as formula (\ref{UP complex formula one in proof}).
Now we prove \[{\rm V}_{x,f}
\int_{\mathbb{S}^n}|\Phi'_f(x)|^2|f(x)|^2d\sigma(x)\geq {\rm COV}^2.\]

By using Cauchy-Schwarz's inequality, we immediately obtain
\begin{eqnarray*}
&&{\rm V}_{x,f}
\int_{\mathbb{S}^n}|\Phi'_f(x)|^2|f(x)|^2d\sigma(x)\\
&=&\int_{\mathbb{S}^n}|x-\tau_f|^2
|f(x)|^2d\sigma(x)
\int_{\mathbb{S}^n}|\Phi'_f(x)|^2|f(x)|^2d\sigma(x)\\
&\geq&[\int_{\mathbb{S}^n}|x-\tau_f|
|f(x)||\Phi'_f(x)||f(x)|d\sigma(x)]^2\\
&=&[\int_{\mathbb{S}^n}|x-\tau_f|
|\Phi'_f(x)||f(x)|^2d\sigma(x)]^2.
\end{eqnarray*}
Hence, 
\begin{align}\label{UP-vector}
    {\rm V}_{x,f}{\rm V}_{\nabla_{\mathbb{S}^n},f}\geq\frac{n^2}{4}|\tau_f|^4+{\rm COV}^2,
\end{align}
where
\begin{align*}
{\rm COV}^2\triangleq [\int_{\mathbb{S}^n}|x-\tau_f|
|\Phi'_f(x)||f(x)|^2d\sigma(x)]^2.
\end{align*}
Note that 
\begin{align*}
    a(f) = (a_0(f),a_1(f),...,a_n(f))^T,
\end{align*}
where $a_k(f)=\sum_{j=1}^m\int_{\mathbb S^n}f_j(x)(\Gamma_kf_j(x)) d\sigma(x).$ Then by Lemma \ref{integration by parts}, we have
\begin{align*}
    a_k(f) = \sum_{j=1}^m \left(n\int_{\mathbb S^n}x_k|f_j(x)|^2d\sigma(x)-\int_{\mathbb S^n}f_j(x)(\Gamma_kf_j(x))d\sigma(x)\right),
\end{align*}
which gives
\begin{align*}
    a_k(f) = \frac{n}{2}\sum_{j=1}^m \int_{\mathbb S^n}x_k |f_j(x)|^2d\sigma(x) =\frac{n}{2}\int_{\mathbb S^n}x_k|f(x)|^2d\sigma(x)
\end{align*}
for $k=0,1,...,n.$
Hence,
\begin{align}\label{af_tauf-vector}
    \left|\int_{\mathbb S^n}(\nabla_{\mathbb S^n}f(x))f^T(x)d\sigma(x)\right|^2 = \frac{n^2}{4} \sum_{k=0}^n\left|\int_{\mathbb S^n}x_k|f(x)|^2d\sigma(x)\right|^2=\frac{n^2}{4}\left|\int_{\mathbb S^n}x|f(x)|^2d\sigma(x)\right|^2.
\end{align}
By (\ref{UP-vector}), we have
\begin{align*}
    &\int_{\mathbb S^n}|\nabla_{\mathbb S^n}f(x)-a(f)f(x)|^2d\sigma(x) \left(1-\left|\int_{\mathbb S^n}x|f(x)|^2d\sigma(x)\right|^2\right)\\
    = & \left(\int_{\mathbb S^n}|\nabla_{\mathbb S^n}f(x)|^2d\sigma(x)-|a(f)|^2\right)\left(1-|\tau_f|^2\right) \\
    \geq &\frac{n^2}{4}|\tau_f|^4+\COV^2.
\end{align*}
Then
\begin{align*}
    &\int_{\mathbb S^n}|\nabla_{\mathbb S^n}f(x)|^2d\sigma(x)\left(1-\left|\int_{\mathbb S^n}x|f(x)|^2d\sigma(x)\right|^2\right)\\
    \geq &\frac{n^2}{4}|\tau_f|^4+\COV^2+|a(f)|^2-|a(f)|^2|\tau_f|^2\\
    = & \frac{n^2}{4}|\tau_f|^2+\COV^2,
\end{align*}
where we have used (\ref{af_tauf-vector}).
\endproof	

\noindent \textbf{Acknowledgements.} W.-X. Mai was supported by NSFC Grant No.11901594, FRG Program of the Macau University of Science and Technology, No. FRG-22-076-MCMS, and the Science and Technology Development Fund, Macau SAR (File no. 0133/2022/A). P. Dang was supported by the Science and Technology Development Fund, Macau SAR (File no. 0006/2019/A1).

\end{document}